\documentclass[letterpaper, 10 pt, conference]{ieeeconf}  %
\IEEEoverridecommandlockouts
\title{\LARGE \bf
Self-sustained oscillations in discrete-time relay feedback systems \thanks{The project was supported by the Israel Science Foundation (grant no. 2406/22), while the second author was a Jane and Larry Sherman Fellow. All the authors acknowledge the EuroTech 2025 Cooperative Scholarship for funding this research.}
}

\author{Kang Tong, Christian Grussler and Michelle S. Chong   \thanks{K. Tong and C. Grussler are with the Stephen B. Klein Faculty of Aerospace Engineering, Technion -- Israel Institute of Technology, Haifa, Israel. Emails: {\tt\small kang.tong@campus.technion.ac.il, cgrussler@technion.ac.il} }
\thanks{M. Chong is with the Department of Mechanical Engineering, Eindhoven University of Technology, Eindhoven, Netherlands. Email: {\tt\small m.s.t.chong@tue.nl}}
\thanks{This paper has been accepted as an invited session paper at the 64th IEEE Conference on Decision and Control (CDC 2025).}
}

\usepackage{graphicx}
\usepackage{enumerate}
\usepackage{mathtools}
\usepackage{xcolor} %

\usepackage{amsmath} %
\usepackage{amssymb}  %
\usepackage{bm}
\usepackage{fixme}
\fxsetup{status=draft}

\usepackage{bbm}
\usepackage{algorithm}
\usepackage{algorithmic}
\usepackage{pgfplots}
\pgfplotsset{compat=newest}

\usepackage{tikz}
\usepackage{tkz-euclide}

\usepackage{enumerate}

\usepackage{mathrsfs}
\definecolor{ao(english)}{rgb}{0.0, 0.5, 0.0}
\usepackage[colorlinks, citecolor = {ao(english)}, linkcolor = {ao(english)}]{hyperref} 
\usepackage[capitalize,nameinlink]{cleveref}
\usepackage{siunitx}
\usepackage[sort,nocompress]{cite}
\usepackage{subcaption}

\usetikzlibrary{arrows,shapes,trees,calc,positioning,patterns,decorations.pathmorphing,decorations.markings,backgrounds}
\usetikzlibrary{matrix}

\usepgfplotslibrary{groupplots}

\newtheorem{thm}{Theorem}
\crefname{thm}{Theorem}{Theorems}

\crefname{prop}{Proposition}{Propositions}

\newtheorem{lem}{Lemma}
\crefname{lem}{Lemma}{Lemmas}

\newtheorem{cor}{Corollary}
\crefname{cor}{Corollary}{Corollaries}

\newtheorem{rem}{Remark}
\crefname{rem}{Remark}{Remark}

\newtheorem{ass}{Assumption}
\crefname{ass}{Assumption}{Assumption}

\crefname{conj}{Conjecture}{Conjectures}

\newtheorem{defn}{Definition}
\crefname{defn}{Definition}{Definitions}

\crefname{prob}{Problem}{Problems}

\crefname{appl}{Application}{Applications}

\newtheorem{exm}{Example}
\crefname{exm}{Example}{Examples}

\crefname{algorithm}{Algorithm}{Algorithms}
\crefname{paper}{Paper}{Papers}
\crefname{figure}{Figure}{Figures}
\crefname{section}{Section}{Sections}

\usepackage{dsfont}
\let\mathbb=\mathds

\newcommand{\sign}{\textnormal{sign}}

\newcommand{\argmin}{\operatornamewithlimits{argmin}}

\colorlet{FigColor1}{blue}
\colorlet{FigColor2}{red}
\colorlet{FigColor3}{ao(english)}
\colorlet{FigColor4}{orange}
\pgfplotsset{every axis plot/.append style={line width=1.5pt}}
\definecolor{bluebell}{rgb}{0.74, 0.83, 0.9}
\definecolor{airforceblue}{rgb}{0.36, 0.54, 0.66}

\crefformat{equation}{\textup{#2(#1)#3}}
\crefrangeformat{equation}{\textup{#3(#1)#4--#5(#2)#6}}
\crefmultiformat{equation}{\textup{#2(#1)#3}}{ and \textup{#2(#1)#3}}
{, \textup{#2(#1)#3}}{, and \textup{#2(#1)#3}}
\crefrangemultiformat{equation}{\textup{#3(#1)#4--#5(#2)#6}}%
{ and \textup{#3(#1)#4--#5(#2)#6}}{, \textup{#3(#1)#4--#5(#2)#6}}{, and \textup{#3(#1)#4--#5(#2)#6}}

\Crefformat{equation}{#2Equation~\textup{(#1)}#3}
\Crefrangeformat{equation}{Equations~\textup{#3(#1)#4--#5(#2)#6}}
\Crefmultiformat{equation}{Equations~\textup{#2(#1)#3}}{ and \textup{#2(#1)#3}}
{, \textup{#2(#1)#3}}{, and \textup{#2(#1)#3}}
\Crefrangemultiformat{equation}{Equations~\textup{#3(#1)#4--#5(#2)#6}}%
{ and \textup{#3(#1)#4--#5(#2)#6}}{, \textup{#3(#1)#4--#5(#2)#6}}{, and \textup{#3(#1)#4--#5(#2)#6}}

\crefdefaultlabelformat{#2\textup{#1}#3}

\begin{document}

\maketitle
\thispagestyle{empty}
\pagestyle{empty}

\begin{abstract}

We study the problem of determining self-sustained oscillations in discrete-time linear time-invariant relay feedback systems. Concretely, we are interested in predicting when such a system admits unimodal oscillations, i.e., when the output has a single-peaked period. Under the assumption that the linear system is stable and has an impulse response that is strictly monotonically decreasing on its infinite support, we take a novel approach in using the framework of total positivity to address our main question. It is shown that unimodal self-oscillations can only exist if the number of positive and negative elements in a period coincides. Based on this result, we derive conditions for the existence of such oscillations, determine bounds on their periods, and address the question of uniqueness.

\end{abstract}

\section{Introduction}

\begin{figure}[t]
	\centering
	\tikzstyle{neu}=[draw, very thick, align = center,circle]
	\tikzstyle{int}=[draw,minimum width=1cm, minimum height=1cm, very thick, align = center]
	\begin{tikzpicture}[>=latex',circle dotted/.style={dash pattern=on .05mm off 1.2mm,
			line cap=round}]
    \node [coordinate, name=input] {};
    \node [int, right of=input, node distance=2cm] (relay) {\begin{tikzpicture}
				\begin{axis}[ticks = none,width = 2.8 cm, axis lines = none]
					\addplot[line width = 1 pt,color = blue] 
                    table{
                    -4 -2
                    0 -2
                    0 2
                    4 2
};
				\end{axis}
		\end{tikzpicture}};
    \node [int, right of=relay, node distance=2cm] (system) {$G(z)$};
    \node [coordinate][right of=system, node distance=2cm] (output) {};
    \node [int, below left of=system, node distance=1.5cm, shift={(0cm, -0.5cm)}] (feedback) {$-1$};
    \draw [-] (feedback) -| node[above] {} node[below] {} (input);
    \draw [->] (input) -- node[above] {$u(t)$} (relay);
    \draw [->] (relay) -- node[above] {} (system);
    \draw [-] (system) -- node[above] [name=y] {$y(t)$} (output);
    \draw [->] (output) |- node[above] {} node[below] {} (feedback);
    \end{tikzpicture}
	
    \caption{Discrete-time relay feedback system with linear system $G(z)$, $z\in\mathbb{C}$. \label{fig:relay_feedback}}
\end{figure}
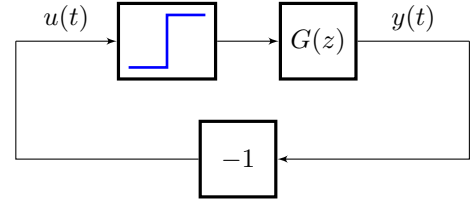

Self-oscillations are common in both natural and engineered systems, ranging from the nervous system \cite{edelstein2005limit, mackey1977oscillation} to robot locomotion \cite{ijspeert2013dynamical} and PID-autotuning with relay feedback \cite{aastrom2004revisiting}. Motivated by the practicality of relay autotuning and the observation that complex oscillatory behavior can occur in relay feedback systems \cite{tsypkin1984relay, holmberg1991relay, bernardo2001self, rabi2020relay}, we focus on self-oscillations in \textit{ideal relay} feedback systems, which consists of a relay function in negative feedback with a linear system as shown in Figure \ref{fig:relay_feedback}. In particular, in predicting the occurrence of self-oscillations for the feedback system in \textit{discrete-time}, which is a challenge and has remained relatively unexplored to date.

So far, a series of tools have been developed for continuous-time relay feedback systems, including graphical tools such as the Hamel locus \cite{le1970general} in time-domain and the Tsypkin locus  \cite{tsypkin1984relay, judd1974graphical, judd1977error} in frequency-domain, which are facilitated by the describing function analysis using sinusoids \cite[Chap. 7.2]{khalil2002nonlinear} and recently, square waves \cite{chaffey2024amplitude}, to provide \textit{approximate} necessary conditions for the existence of self-oscillations. Since then, other approaches to provide exact conditions were developed, such as Lyapunov-based analysis in \cite{gonccalves2001global}; fixed points analysis of Poincaré map - a mapping from the state at one mode to a different mode \cite{varigonda2001dynamics}; and characterizing the linear system itself \cite{megretski1996global, johansson1999fast}.

The literature described thus far has been for continuous-time relay feedback systems. To the best of our knowledge, self-oscillations for discrete-time relay feedback systems have not been studied. In fact, current approaches for Lur'e systems in discrete-time -- a linear system in feedback with a static nonlinearity, proposed in \cite{rasvan1998self, das2020analysis, paredes2024self}, are not applicable to the relay function, as it is neither piecewise continuous, nor does it belong to the conic sector. This motivates us to develop a new approach in characterizing self-oscillations for discrete-time relay feedback system (DT-RFS), depicted in \cref{fig:relay_feedback}.

In this paper, we introduce a novel approach using the total positivity framework of Karlin \cite{karlin1968total} and one of its recent extensions \cite{grussler2025tractable} to rigorously describe and classify oscillatory patterns via the notion of cyclic variation (number of sign changes per period). This enables us to derive broader conclusions about a series of self-oscillations that share similar variation patterns.
In particular, note that the unimodality of self-oscillations -- as observed in \cite{rabi2020relay, gonccalves2001global, varigonda2001dynamics} -- corresponds to a cyclic variation of the incremental sequence equal to two.

We first examine the characteristics of unimodal self-oscillations when the impulse response of the discrete-time linear time-invariant (DT-LTI) system in the DT-RFS is strictly monotonically decreasing on its infinite support. It is shown that unimodal self-oscillations can only exist if the number of positive and negative elements in a period coincide. Based on this result, we derive conditions under which the DT-RFS does not exhibit oscillations when the system is of relative degree zero.

More generally, we demonstrate that the maximum oscillation period is bounded by the delay time, which is derived from the impulse response. Given the maximum period, one can determine the existence of other oscillation periods. Furthermore, if the impulse response is convex on its support, we can obtain a tighter bound for the maximum period. These conclusions are verified by several examples.

The remainder of the paper is organized as follows. After some preliminaries in \cref{sec:prelim}, we state our problem in \cref{sec:prob} and provide the characteristics of self-oscillations in \cref{sec:charac_var} based on the theory of variation. Subsequently, in \cref{sec:main_res}, we derive our main results on the characterization of unimodal self-oscillation with strictly monotonically decreasing impulse response and illustrate them with examples. Conclusions are presented in \cref{sec:conclusion}. 
All proofs are provided in the appendix.

\section{Preliminaries} \label{sec:prelim}
In the following, we introduce several notations and concepts needed for our subsequent derivations and discussions. 
\subsection{Notations}
\subsubsection{Sets} We write $\mathbb{C}$ for the set of complex numbers, $\mathbb{Z}$ for the set of integers and $\mathbb{R}$ for the set of reals, with $\mathbb{Z}_{\geq0}$ ($\mathbb{Z}_{>0}$) and $\mathbb{R}_{\geq0}$ ($\mathbb{R}_{>0}$) standing for the respective subsets of nonnegative (positive) elements. For $k, l \in \mathbb{Z}$ with $k \leq l$, we write $(k:l) := \{ 
k, k+1, \cdots, l \}$.
\subsubsection{Sequences} For a sequence $x : \mathbb{Z} \to \mathbb{R}$, with $\sup_{i \in \mathbb{Z}} |x(i)| < \infty$ or $\sum_{i \in \mathbb{Z}} |x(i)| < \infty$, we write $x \in \ell_{\infty}$ or $x \in \ell_{1}$, respectively. We use $\Delta x(t):= x(t+1)-x(t)$ to denote the \emph{forward difference operator} for a sequence $x$ and say that $x$ is \emph{convex} if $\Delta x(t) \geq \Delta x(t-1)$ for all $t \in \mathds{Z}$. Further, $x$ is called \emph{(strictly) unimodal}, if there exists an $m \in \mathbb{Z}$ for which $x(i)$ is (strictly) monotonically increasing for all $i \leq m$ and (strictly) monotonically decreasing for $t \geq m$. For a slice of the sequence $x$, we define the vector $x(k:l):= \begin{bmatrix}
    x(k) & x(k+1) & \dots & x(l)
\end{bmatrix}^\top$ for $k \leq l$. 

If there exists a $P \in \mathbb{Z}_{>0}$ such that $x(i) = x(i+P)$ for all $i$, then $x$ is called \emph{$P$-periodic.} The set of all bounded $P$-periodic sequences is denoted by $\ell_\infty(P)$.
Additionally, $P$-periodic $x$ is called \emph{(strictly) periodically unimodal} if there exist $k_1, k_2 \in \{0,\dots, P-1\}$ for which $x(i+k_1)$ is (strictly) monotonically increasing for all $0 \leq i \leq k_2$ and (strictly) monotonically decreasing for $k_2 \leq i \leq P-1$. 
Furthermore, the $P$-truncation vector of $x \in \ell_\infty(P)$ over a single period is denoted by
\begin{align*}
    x^P := \begin{bmatrix}
        x(0) \\ \vdots \\ x(P-1)
    \end{bmatrix}.
\end{align*}
\subsubsection{Matrices} For matrices $M = (m_{ij}) \in \mathds{R}^{n \times m}$, we say that $M$ is
\emph{nonnegative},
$M \geq 0$ 
or $M \in \mathds{R}^{n \times m}_{\geq 0}$
if all elements $m_{ij} \in \mathbb{R}_{\geq 0}$ -- corresponding notations are also used for 
$M$ is \emph{positive}, $M>0$ or $M \in \mathds{R}^{n \times m}_{> 0}$ if all elements $m_{ij} \in \mathbb{R}_{> 0}$.
Let $I_n$ be the identity matrix in $\mathds{R}^{n \times n}$ and $\boldsymbol{1}_n$ be the vector of all ones in $\mathbb{R}^n$.  For a vector $v = (v_i) \in \mathbb{R}^n$, the \emph{cyclic forward difference operator} $\Delta_c$ is defined as 
\begin{align*}
    \Delta_c v := \Delta_c \begin{bmatrix}
        v_1 \\ v_2 \\ \vdots \\ v_n
    \end{bmatrix} = \begin{bmatrix}
        \Delta_c v_1 \\ \Delta_c v_2 \\ \vdots \\ \Delta_c v_n
    \end{bmatrix} = \begin{bmatrix}
        v_2 - v_1 \\ v_3 - v_2 \\ \vdots \\ v_1 - v_n
    \end{bmatrix}
\end{align*}
and the \emph{circulant matrix} generated from $v$ is denoted by
\begin{align} \label{eq:circulant_mat}
    H_v := \begin{bmatrix}
        v_1 & v_n & \cdots & v_2 \\
        v_2 & v_1 & \cdots & v_3 \\
        \vdots & \vdots &  & \vdots \\
        v_n & v_{n-1} & \cdots & v_1
    \end{bmatrix} \in \mathbb{R}^{n \times n}.
\end{align}
Finally, we denote the \emph{cyclic back shift matrix of order $n$} by
\begin{align}
    Q_n := \begin{bmatrix} \label{eq:cyclic_back_shift_mat}
        0 & 1 \\
        I_{n-1} & 0  
    \end{bmatrix} \in \mathbb{R}^{n \times n}.
\end{align}

\subsubsection{Functions} The \emph{indicator function} of set $\mathcal{A} \subset \mathbb{R}$ is denoted by $\mathbb{1}_{\mathcal{A}}(t)=1$ for $t\in \mathcal{A}$ and $\mathbb{1}_{\mathcal{A}}(t)=0$ for $t \notin \mathcal{A}$.
In terms of the indicator function $\mathbb{1}_{\mathcal{A}}$, the \emph{unit pulse function} is denoted by $\delta(t):=\mathbb{1}_{\{0\}} (t)$ and the \emph{sign function} or \emph{ideal relay function} is denoted by
\begin{align*}
    \sign(t) := - \mathbb{1}_{\mathbb{R}_{< 0}}(t) + \mathbb{1}_{\mathbb{R}_{> 0}}(t).
\end{align*} 
The number of positive, negative, and zero elements in $v \in \mathbb{R}^n$ are defined by
$P_p(v) := \sum_{i=1}^{n} \mathbb{1}_{\mathbb{R}_{> 0}} (v_i)$, $P_n(v) := \sum_{i=1}^{n} \mathbb{1}_{\mathbb{R}_{< 0}} (v_i)$ and $P_z(v) := \sum_{i=1}^{n} \mathbb{1}_{\{0\}} (v_i)$, respectively.
The \emph{ceiling function} $\lceil x \rceil$ maps $x \in \mathds{R}$ to the smallest integer greater than or equal to $x$.

\subsection{Linear Time-Invariant Systems} 

We consider finite-dimensional {causal linear discrete-time invariant (LTI) systems} with (scalar) input $u \in \ell_\infty$ and (scalar) outputs $y \in \ell_\infty$. The \emph{impulse response} $g \in \ell_1$ is the output corresponding to the input $u(t) = \delta(t)$ and the \emph{transfer function} of the system is given by
\begin{equation} \label{eq:def_trans_fun}
G(z) = \sum_{t=0}^\infty g(t)z^{-t}, \quad z \in \mathbb{C}.
\end{equation}
The \emph{convolution operator} $\mathcal{C}_g: \ell_\infty \to \ell_\infty$ with respect to $G$ is defined by 
\begin{align}
   \mathcal{C}_g(u)(t) := (g \ast u) (t) := \sum_{\tau = -\infty}^{\infty} g(t-\tau) u(\tau), \ t \in \mathbb{Z}.
\end{align}
For $u \in \ell_\infty(P)$, it holds then that $y \in \ell_\infty(P)$ is such that
\begin{equation} \label{eq:original_conv_cal}
    \begin{aligned}
    y(t) = \mathcal{C}_g(u)(t) & = \sum_{\tau = -\infty}^{\infty} g(t-\tau) u(\tau)  \\
    & = \sum_{j=0}^{P-1} \left( u(j) \sum_{m = -\infty}^{\infty} g(t+mP-j) \right)  \\
    & = \sum_{j=0}^{P-1} u(j) \overline{g}(t-j), 
\end{aligned}
\end{equation}
with $\overline{g}(t) := \sum_{i=-\infty}^{\infty} g(t + iP)$
denoting the so-called \emph{periodic summation} of $g$. Equivalently, 
\begin{align}
    y^P =\begin{bmatrix} \label{eq:conv_to_mat_mul}
        \mathcal{C}_g(u)(0) \\ \vdots \\ \mathcal{C}_g(u)(P-1)
    \end{bmatrix} = H_{\overline{g}^P} u^P = H_{u^P} \overline{g}^P.
\end{align}

\section{Problem Statement} \label{sec:prob}
In this work, we study the existence of sustained unimodal oscillations in discrete-time relay feedback systems (see~\cref{fig:relay_feedback}), which is comprised of a causal linear-time invariant system $G(z)$ in negative feedback with the ideal relay function. Using our defined operator notations, this feedback interconnection reads as
\begin{align} \label{eq:closed_loop_sys}
    u(t) = -\mathcal{C}_g (\sign(u))(t), \quad t \in \mathbb{Z}.
\end{align}
Therefore, unimodal self-sustained oscillations can be identified with periodically unimodal fixed-points $u \in \ell_\infty(P)$ satisfying \cref{eq:closed_loop_sys}.
In particular, we make the following assumptions on the impulse response $g$:
\begin{ass} \label{ass:mono_dec_g}
    The impulse response $g$ to $G(z)$ satisfies $g \in \ell_{1}$, has a connected support $\{t_1 + i\}_{i \in \mathbb{Z}_{\geq 0}}$ for some $t_1 \in \mathbb{Z}_{\geq 0}$ and is strictly monotonically decreasing on its support.  
\end{ass}
This assumption will be crucial to ensure that the system leaves the set of periodically unimodal signals invariant. Note that since $g \in \ell_1$, it also follows that $g$ is strictly positive on its support.

\section{Self-oscillations via variation} \label{sec:charac_var}
A main tool for our studies is the framework of total positivity \cite{karlin1968total,grussler2025tractable}, which characterizes convolution operators $\mathcal{C}_g$ that map periodically unimodal inputs to periodically unimodal outputs via the notion of variation. In this section, we define the required notions and review relevant preliminary results of that framework. 
\begin{defn} [Variation] 
    For a vector $v \in \mathbb{R}^n$, we define $S^-(v)$ as the number of sign changes in $v$, i.e., 
    \begin{align*}
        S^-(v) := \sum_{i=1}^{m-1} \mathbb{1}_{\mathbb{R}_{<0}} (\tilde{v}_{i} \tilde{v}_{i+1}), \quad S^-(0) = -1,
    \end{align*}
    where $\tilde{v} \in \mathbb{R}^{m}$ is the vector resulting from deleting all zeros in $v$. Further, we define $S^{+}(v) := S^{-}(\bar{v})$, where $\bar{v} \in \mathbb{R}^n$ results from replacing zero elements in $v$ by any real elements that maximize $S^{-}(\bar{v})$.
\end{defn}
Clearly, $S^{-}(v) \leq S^{+}(v)$, but equality does not necessarily hold. For example, if $v = [1, 0, 3]^\top$, then $S^{-}(v)=0$, but if we choose $\bar{v} = [1, a, 3]^\top$ with $a < 0$, we obtain $S^{+}(v) = S^{-}(\bar{v}) = 2$.
We also need the notion of so-called \emph{cyclic variation} for $v \in \mathbb{R}^n$, which is defined as follows.
\begin{align*}
    S_c^-(v) := \sup_{i \in (1:n)} S^-([v_i, \dots, v_n, v_1, \dots, v_i]), \\
    S_c^+(v) := \sup_{i \in (1:n)} S^+([v_i, \dots, v_n, v_1, \dots, v_i]).
\end{align*}
We are ready, now, to precisely define self-sustained oscillation via  cyclic variation:
\begin{defn} [Self-oscillation] \label{def:osci_orig}
    The sequence $u \in \ell_{\infty}(P)$ is called a self-sustained oscillation if and only if $u$ satisfies \cref{eq:closed_loop_sys} and $S_c^-(\Delta_c u^P) \geq 2$.
\end{defn}
For the remainder of the paper, we will also use the abbreviated term ``self-oscillation". The condition $S_c^-(\Delta_c u^P) \geq 2$ ensures that the sequence $u$ is not constant. This is analogous to the notion of \emph{limit cycles} \cite{khalil2002nonlinear}, where equilibrium points, i.e., constant limit cycles, are excluded. 

Periodic unimodality of $u \in \ell_\infty(P) \setminus \{0\}$ can be equivalently characterized by $S_c^-(\Delta_c u^P) = 2$ (see~\cite[Lem.~3]{grussler2025tractable}). In our analysis, we require the following slightly stronger form of periodic unimodality:
\begin{ass} \label{ass:osci_constr_1}
    The sequence $u \in \ell_{\infty}(P)$ satisfies 
    \begin{align} 
    S_c^-(\Delta_c u^P) = 2  \text{ and } S_c^+(u^P) = 2.
\end{align}
\end{ass}
The following lemma provides a way of verifying whether the circulant matrix $H_v$ satisfies $S_c^-(\Delta_c (H_v w)) \leq 2$ for any $w$ fulfilling $S_c^-(\Delta_c w) = 2$. This is important, as it allows us to check via \cref{eq:conv_to_mat_mul} whether the
\emph{loop gain} 
\begin{align} \label{eq:dt_rfs_without_delay}
    o(t) = -\mathcal{C}_g(\sign(u))(t)
\end{align} 
maps periodically unimodal inputs to periodically unimodal outputs. 
\begin{lem} \label{lem:sign_PM}
    Let $v \in \mathbb{R}^n$ and $S_c^+(v) \leq 2$. Then, for any $w \in \mathbb{R}^n$ satisfying $S_c^-(\Delta_c w) \leq 2$, it holds that
    \begin{align*}
        S_c^-(\Delta_c H_{\sign(v)} w) \leq 2.
    \end{align*}
\end{lem}

\section{Main Results} \label{sec:main_res}
In this section, we present our main results on self-oscillations in DT-RFS. We begin by deriving necessary conditions for the existence of such oscillations. Based on this, it is subsequently shown that there cannot be any self-oscillations if $g(0) > 0$. We then continue discussing the cases, where $g(0) = 0$ for which we provide lower and upper bounds for the period of self-oscillations and address the question of its uniqueness.

\subsection{Necessary conditions for self-sustained oscillations}
To analyze the self-sustained oscillation in the closed-loop system, we begin by showing that if $g$ satisfies \cref{ass:mono_dec_g}, then the loop-gain \cref{eq:dt_rfs_without_delay} is invariant to signals in \cref{ass:osci_constr_1}.
\begin{lem}  \label{lem:Mickey_pres}
    Let $P > 1$, $g$ under \cref{ass:mono_dec_g}, $u \in \ell_{\infty}(P)$ under \cref{ass:osci_constr_1} and let the one period vector of \eqref{eq:dt_rfs_without_delay} be given by 
    \begin{align}  \label{eq:open_loop_rfs}
        o^P = -H_{\overline{g}^P}\sign(u^P).
    \end{align}
    Then, 
    \begin{align*}
        S^-_c(\Delta_c o^P) \leq 2 \text{ and } S^{-}_c(o^P) \leq 2.
    \end{align*}
    If additionally, $P_p(u^P) = P_n(u^P)$, then $S^{+}_c(o^P) \leq 2$, i.e., $o^P$ satisfies \cref{ass:osci_constr_1}.
\end{lem}

Our analysis of the open-loop system properties yields the first main result: a necessary condition for the open-loop system to admit a fixed-point. When the sequence is a fixed-point for the open-loop system, the output and input sequences must share the same number of positive, negative, and zero elements, denoted by $P_p$, $P_n$, and $P_z$, in one period.

\begin{thm} \label{thm:osci_pres_open_loop}
    Let $P>1$, consider $g$ under \cref{ass:mono_dec_g}, $u \in \ell_{\infty}(P)$ under \cref{ass:osci_constr_1} and one period output of open-loop system \eqref{eq:open_loop_rfs}. If $u^P$ and $o^P$ satisfy 
    \begin{align*}
        P_p(u^P) = P_p(o^P) \text{ and } P_n(u^P) = P_n(o^P),
    \end{align*}
    then 
    \begin{align*}
        P_p(u^P) = P_n(u^P).
    \end{align*}
    If additionally $P_z(u^P) = 0$, then 
    \begin{align*}
        P = 2P_p(u^P) = 2P_p(o^P).
    \end{align*} 
\end{thm}

It is important to note that from \cref{lem:Mickey_pres} and \cref{thm:osci_pres_open_loop}, it can be concluded that unimodal self-oscillations satisfying \cref{ass:osci_constr_1} can only exist if the numbers of positive and negative elements in a period are equal.

\subsection{Absence of self-oscillations}
Next, we state our second main result on the absence of unimodal self-oscillation in the case of $g(0)>0$.
\begin{thm} \label{thm:non_oscillation}
    Let $g \in \ell_1$ satisfy $g(0)>0$ and \cref{ass:mono_dec_g}. Then, there does not exist any $u$ satisfying \cref{ass:osci_constr_1} such that $u = -\mathcal{C}_g(\sign(u))$.
\end{thm}
An example for systems satisfying \cref{thm:non_oscillation} is the parallel interconnection of first-order lags,
\begin{align*}
    G(z) = \sum_{i=1}^{n}\frac{k_i z}{z - p_i}
\end{align*}
with $p_i \in (0, 1)$ and $k_i > 0$, where $g(0) = \sum_{i=1}^{n} k_i >0$. In this case, our result closely resembles the absence of self-oscillations known from continuous-time passive systems \cite[Chap. 6]{khalil2002nonlinear}.

\subsection{Self-oscillations under time-delays}

\begin{figure}[t]
	\centering
	\tikzstyle{neu}=[draw, very thick, align = center,circle]
	\tikzstyle{int}=[draw,minimum width=1cm, minimum height=1cm, very thick, align = center]
	\begin{tikzpicture}[>=latex',circle dotted/.style={dash pattern=on .05mm off 1.2mm,
			line cap=round}]
    \node [coordinate, right of=input] (sum) {};
    \node [int, right of=input, node distance=3cm] (relay) {\begin{tikzpicture}
				\begin{axis}[ticks = none,width = 2.8 cm, axis lines = none]
					\addplot[line width = 1 pt,color = blue] table{
                    -4 -2
                    0 -2
                    0 2
                    4 2
                    };
				\end{axis}
		\end{tikzpicture}};
    \node [int, right of=relay, node distance=2cm] (system) {$G_0(z)$};
    \node [coordinate][right of=system, node distance=2cm] (output) {};
    \node [int, below of=system, node distance=1.5cm] (feedback) {$z^{-P_d}$};
    \node [int, below of=relay, node distance=1.5cm] (delay) {$-1$};
    \draw [->] (sum) -- node[above] {$u(t)$} (relay);
    \draw [->] (relay) -- node[above] {} (system);
    \draw [-] (system) -- node[above] [name=y] {$\tilde{y}(t)$} (output);
    \draw [->] (output) |- (feedback);
    \draw [->] (feedback) -- (delay);
    \draw [-] (delay) -| node[pos=0.9, left] {} node [near end, left] {} (sum);
    \end{tikzpicture}
	
    \caption{A DT Relay Feedback System with delay module. \label{fig:dt_rfs_with_delay}}
\end{figure}
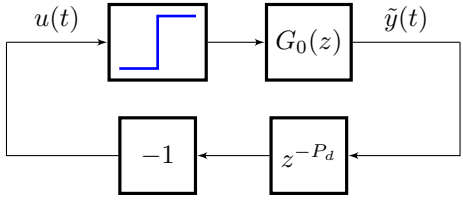

Next, we will show that self-oscillations exist if $g$ satisfies \cref{ass:mono_dec_g} with $g(0) = 0$.
In this case, we can represent $G(z) = z^{-P_d} G_0(z)$, where $P_d \geq 0$ and $G_0(z)$ is causal with $g_0(0) > 0$.
This conversion is helpful for analyzing the contributions of $g_0$ and the time-delay $P_d$ independently, where the closed-loop system depicted in \cref{fig:dt_rfs_with_delay} reads as 
\begin{align} \label{eq:dt_rfs_with_delay}
    u(t) & = -\mathcal{C}_g(\sign(u))(t) \nonumber \\
    & = -\mathcal{C}_{g_0}(\sign(u))(t-P_d).
\end{align} 
or equivalently, 
\begin{align} \label{eq:osci_delay_mat_form}
    u^P & = \begin{bmatrix}
        -\mathcal{C}_{g_0}(\sign(u))(-P_d) \\ \vdots \\ -\mathcal{C}_{g_0}(\sign(u))(P-P_d-1)
    \end{bmatrix} \nonumber \\
    & = -Q_{P}^{P_d} H_{\overline{g}_0^P} \sign(u^P).
\end{align}
We are now ready to present our third main result: the relationship between $P_d$ and self-oscillations with period $P \geq P_d > 0$.

\begin{thm}\label{thm:bound_max_osc}
Consider system \eqref{eq:dt_rfs_with_delay} with $P_d \geq 1$ and $g_0$ satisfying \cref{ass:mono_dec_g}. Then, for any self-oscillation $u \in \ell_{\infty}(P)$ satisfying \cref{ass:osci_constr_1} and $P \geq P_d$, it must hold that
\begin{align*}
    2 P_d \leq P \leq 2 (P_d+P_s),
\end{align*}
where 
\begin{align*}
    P_s = \argmin_{t} \left(\sum_{k=0}^{t-1} g_0(k)- \sum_{k=t}^{\infty} g_0(k) >0 \right).
\end{align*}
In particular, $P = 2P_d$ if and only if 
\begin{align*}
    \mathcal{C}_{g_0} (\sign{(u)}) (0) > 0,
\end{align*}
and 
$\sign(u^{P}) = \sign(u^{2P_d}) = [\boldsymbol{1}_{P_d}^\top, -\boldsymbol{1}_{P_d}^\top]^\top$.
\end{thm}
\begin{rem}
    A direct conclusion for this theorem is that no self-oscillation whose period is between $P_d$ and $2P_d$.
\end{rem}
Next, we show how \cref{thm:bound_max_osc} can be used to address the question of having more than one self-oscillation at the same time.  

\begin{cor} \label{cor:max_per_to_other}
    Consider system \eqref{eq:dt_rfs_with_delay} where $P_d \geq 1$ and $g_0$ satisfies \cref{ass:mono_dec_g}. If there exists a self-oscillation with period $2P_d$, then there also exist self-oscillations $u \in \ell_{\infty}(P)$ with sign pattern $\sign(u^{P}) = [\boldsymbol{1}_{\frac{P}{2}}^\top, -\boldsymbol{1}_{\frac{P}{2}}^\top]^\top$ for all (even) periods
    \begin{align*}
        P = \frac{2P_d}{2n+1} 
    \end{align*}
    whenever $P \in \mathbb{Z}_{>0}$ for $n \in \mathbb{Z}_{\geq 0}$.
\end{cor}

\begin{exm}
    Consider system \eqref{eq:dt_rfs_with_delay} and let 
    \begin{align*}
        G_0(z)= \frac{z}{z-0.1}.
    \end{align*}
    We plot the relationship between $P$ and $P_d$ by verifying \eqref{eq:osci_delay_mat_form} under the sign pattern $\sign(u^{P}) = [\boldsymbol{1}_{\frac{P}{2}}^\top, -\boldsymbol{1}_{\frac{P}{2}}^\top]^\top$. In \cref{fig:plot_osci_period}, the red points represent the existing oscillation of period $2P_d$ and the blue points represent other existing periods. 
    According to the results in \cref{fig:plot_osci_period}, this system exhibits self-oscillations with periods of 2, 6, and 18 for $P_d = 9$, as illustrated in \cref{fig:osci_period_P_d_9} and verified through simulation.

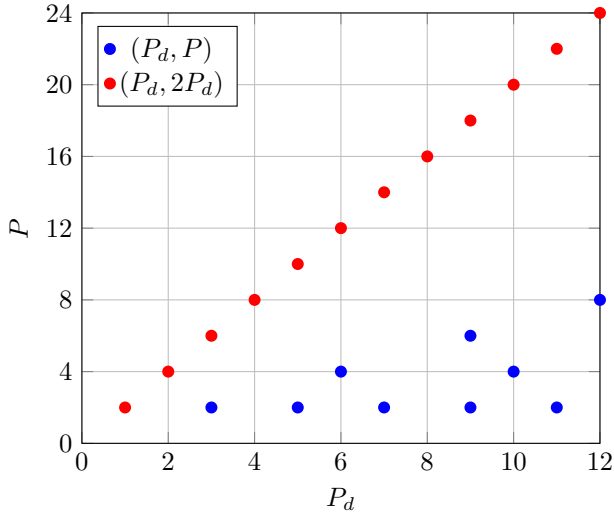
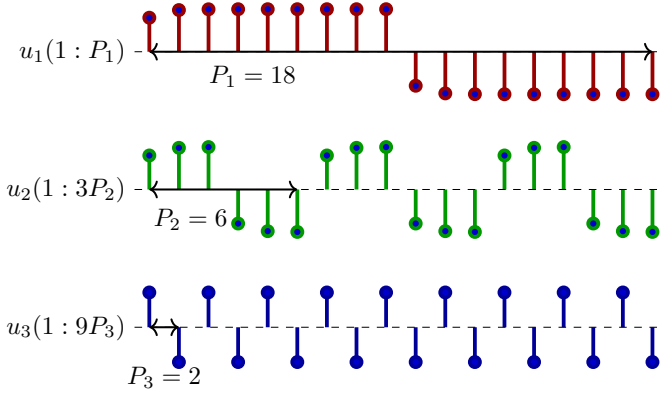
\begin{figure}[t]  
    \begin{subfigure}[b]{\linewidth}
        \centering
        \begin{tikzpicture}
            \begin{axis}[
                xlabel=$P_d$,
                ylabel=$P$,
                xmin=0, xmax=12,
                ymin=0, ymax=24,
                xtick={0,2,4,6,8,10,12},
                ytick={0,4,8,12,16,20,24},
                grid=both, %
                legend pos= north west, %
            ]
            \addplot[
                only marks, %
                color=blue, %
                mark=*, %
                mark size=1.5pt, %
            ] table { %
                x   y
                3   2
                5   2
                6   4
                7   2
                9   2
                9   6
                10  4
                11  2
                12  8
            };
            \addlegendentry{$(P_d,P)$}; %
    
            \addplot[
                only marks, %
                color=red, %
                mark=*, %
                mark size=1.5pt, %
            ] table { %
                x   y
                1   2
                2   4
                3   6
                4   8
                5   10
                6   12
                7   14
                8   16
                9   18
                10  20
                11  22
                12  24
            };
            \addlegendentry{$(P_d,2P_d)$}; %
            \end{axis}
        \end{tikzpicture}
        \caption{Relationship between $P$ and $P_d$.}
        \label{fig:plot_osci_period}
    \end{subfigure}
    \par
    \vspace{1em}
    \begin{subfigure}[b]{\linewidth}
        \centering
    \begin{tikzpicture}[scale=0.95]
        \centering
        \begin{axis}[
            name=plot1,
            width=9cm, height=3.5cm,
            axis lines=none,         %
            ticks=none,              %
            xlabel={}, ylabel={},    %
            title={},                %
            axis on top=true,        %
            clip=false,              %
            xmin=0, xmax=18,
            ymin=-1.8, ymax=1.8,
        ]
        \addplot+[ycomb, color = red!60!black] 
        table {osci_P_9.txt};
        \draw[dashed] (0,0) -- (18,0);
        \node[anchor=east] at (0,0) {$u_1(1:P_1)$};
        \draw[<->, thick] (0.5,0) node[below right =2pt and 20 pt] {$P_1=18$} -- (17.5,0);
        \end{axis}

        \begin{axis}[
            at={(plot1.south)},
            name=plot2,
            anchor=north,
            width=9cm, height=3.5cm,
            axis lines=none,         %
            ticks=none,              %
            xlabel={}, ylabel={},    %
            title={},                %
            axis on top=true,        %
            clip=false,              %
            xmin=0, xmax=18,
            ymin=-1.8, ymax=1.8,
        ]
          \addplot+[ycomb, color = green!60!black] 
            table {osci_P_3.txt};
                
            \draw[dashed] (0,0) -- (18,0);
            \node[anchor=east] at (0,0) {$u_2(1: 3P_2)$};
            \draw[<->, thick] (0.5,0)node[below right =4pt and -2pt] {$P_2=6$} -- (5.5,0);
      \end{axis}

      \begin{axis}[
            at={(plot2.south)},
            name=plot3,
            anchor=north,
            width=9cm, height=3.5cm,
            axis lines=none,         %
            ticks=none,              %
            xlabel={}, ylabel={},    %
            title={},                %
            axis on top=true,        %
            clip=false,              %
            xmin=0, xmax=18,
            ymin=-1.8, ymax=1.8,
        ]
          \addplot+[ycomb, color = blue!60!black] 
            table {osci_P_1.txt};
            \draw[dashed] (0,0) -- (18,0);
            \node[anchor=east] at (0,0) {$u_3(1: 9P_3)$};
            \draw[<->, thick] (0.5,0) -- (1.5,0)
          node[midway, below=12pt] {$P_3=2$};
      \end{axis}

    \end{tikzpicture}
    \caption{Three different self-oscillations when $P_d = 9$.}
    \label{fig:osci_period_P_d_9}
    \end{subfigure}
    
    \caption{Self-oscillations for Example 1.}
\end{figure}

\end{exm}
From this example, we see that in many cases, the oscillation is not unique. The bounds in \cref{thm:bound_max_osc} can be further tightened if we additionally assume that $g_0$ is convex on its support, e.g., $g_0(t) = a^t$ with $a \in (0,1)$ for $t \in \mathbb{Z}_{\geq 0}$.

\begin{cor} \label{cor:convex_bound}
    Consider system \eqref{eq:dt_rfs_with_delay} where $P_d > 1$ and $g_0$ satisfies \cref{ass:mono_dec_g} and $g_0$ being convex on its support. 
    Then, for any self-oscillation $u \in \ell_{\infty}(P)$ satisfying \cref{ass:osci_constr_1} with $P \geq P_d$, it must hold that 
    \begin{align*}
        2 P_d \leq P \leq 4 P_d +2.
    \end{align*}
\end{cor}

The following example verifies that the convex impulse response function $g_0$ satisfies the tighter bound given in \cref{cor:convex_bound}.

\begin{exm}
    Consider system \eqref{eq:dt_rfs_with_delay} with two different linear systems 
    \begin{align*}
        G_{0,1}(z)= \frac{z}{z-0.1} \text{ and } G_{0,2}(z)= \frac{z}{z-0.9},
    \end{align*} 
    with corresponding convex impulse responses
    \begin{align*}
        g_{0,1}(t) = 0.1^{t} \text{ and } g_{0,2}(t) = 0.9^{t}, \text{ for } t \in \mathbb{Z}_{\geq 0}.
    \end{align*}
    For $P \geq P_d$, we plot the relationship between $P$ and $P_d$ by verifying \eqref{eq:osci_delay_mat_form} under the sign pattern $\sign(u^{P}) = [\boldsymbol{1}_{\frac{P}{2}}^\top, -\boldsymbol{1}_{\frac{P}{2}}^\top]^\top$. In \cref{fig:plot_osci_max}, the blue points represent the largest periods of $G_{0,1}(z)$ and the red points represent the largest periods of $G_{0,2}(z)$. The two dashed lines represent the lower bound $2P_d$ and the upper bound $4P_d+2$, respectively. This figure illustrates the conclusion in \cref{cor:convex_bound}.

\begin{figure}[t]
    \centering   
    \begin{tikzpicture}
        \begin{axis}[
            xlabel=$P_d$,
            ylabel=$P$,
            xmin=0, xmax=8,
            ymin=0, ymax=32,
            xtick={0,1,2,3,4,5,6,7,8},
            ytick={0,4,8,12,16,20,24,28,32},
            grid=both, %
            legend pos= north west, %
        ]
        \addplot[
            only marks, %
            color=blue, %
            mark=*, %
            mark size=1.5pt, %
        ] table { %
            x   y
            1   2
            2   4
            3   6
            4   8
            5   10
            6   12
            7   14
            8   16
        };
        \addlegendentry{$G_{0,1}(z) = \frac{z}{z-0.1}$}; %

        \addplot[
            only marks, %
            color=red, %
            mark=*, %
            mark size=1.5pt, %
        ] table { %
            x   y
            1   2
            2   6
            3   10
            4   14
            5   18
            6   22
            7   26
            8   30
        };
        \addlegendentry{$G_{0,2}(z) = \frac{z}{z-0.9}$}; %

        \addplot[magenta, dashed] coordinates{(0,2) (7.5,32)};
        \addlegendentry{$P=4P_d+2$}; %

        \addplot[cyan, dashed] coordinates{(0,0) (8,16)};
        \addlegendentry{$P=2P_d$}; %

        \end{axis}
    \end{tikzpicture}
    \caption{Relationship between $P$ and $P_d$ for different convex impulse responses under $P\geq P_d$.}
    \label{fig:plot_osci_max}
\end{figure}
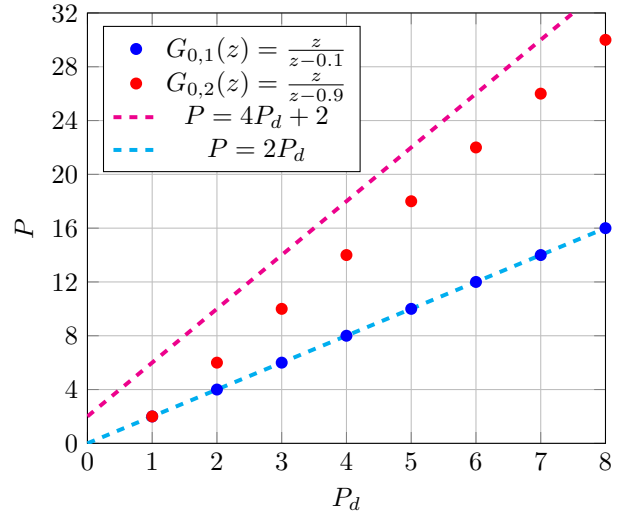

\end{exm}

\section{Conclusion} \label{sec:conclusion}
In this work, we analyze the existence of unimodal self-oscillation in DT-RFS under the assumption of a strictly monotonically decreasing impulse response. Leveraging the framework of total positivity, we show that if the DT-RFS exhibits a unimodal self-oscillation, the number of positive and negative elements in one period must be equal. Based on this result, we provide conditions for the absence and existence of such oscillations as well as bounds for their periods. 
Ultimately, we believe that our analysis in the discrete-time domain will play a crucial role in understanding the continuous-time domain through the process of discretization. In the future, it would also be interesting to explore oscillations with higher variation and to investigate their relationship with the pattern of the impulse response.

\appendix
\subsection{Auxiliary Results}

\begin{lem}[\cite{grussler2025tractable}] \label{lem:VB_2}
    Consider $v \in \mathbb{R}^n$ with $n > 3$ satisfying the following conditions:
    \begin{enumerate}
        \item $S_c^-(\Delta_c v) \leq 2$.
        \item For any $t \in (1:n)$, we have $(\Delta_c v_t)^2 \geq \Delta_c v_{t-1} \Delta_c v_{t+1}$, where the indices are taken modulo $n$; that is, $\Delta_c v_{t+kn} = \Delta_c v_{t}$ for all $k \in \mathbb{Z}$. \label{item:PMP_ineq}
    \end{enumerate}
    Then, $S_c^-(\Delta_c H_{v} w) \leq 2$ for any $w \in \mathbb{R}^n$ that satisfies $S_c^-(\Delta_c w) \leq 2$.
\end{lem}
Note that if $u$ satisfies \cref{ass:osci_constr_1}, then we are only required to consider the following four main potential forms of $\sign(u^P)$: 
\begin{subequations} \label{eq:four_cases}
    \begin{align}
        & [1,\cdots,1,0,-1,\cdots,-1,0]^\top, \label{eq:case_a} \\
        & [1,\cdots,1,-1,\cdots,-1,0]^\top, \label{eq:case_b} \\
        & [1,\cdots,1,0,-1,\cdots,-1]^\top, \label{eq:case_c} \\ 
        & [1,\cdots,1,-1,\cdots,-1]^\top. \label{eq:case_d}
    \end{align}
\end{subequations}
All remaining cases can be generated by those four forms via multiplication with the cyclic back-shift matrix $Q_P^k$ for some $k \in \mathbb{Z}$. Thus, in the following proofs we only need to discuss the four cases stated in \eqref{eq:four_cases}.

\subsection{Proof of \cref{lem:sign_PM}}

    We first consider the cases $n=2$ and $n=3$.
    For these dimensions, every vector $v_a \in \mathbb{R}^n$ satisfies $S_c^-(v_a) \leq 2$. Hence, since $\Delta_c H_{\sign(v)} w \in \mathbb{R}^n$, we have $S_c^-(\Delta_c H_{\sign(v)} w \in \mathbb{R}^n) \leq 2$ for any $v \in \mathbb{R}^n$ and $w \in \mathbb{R}^n$.

    Next, we consider the cases $n>3$. Since $v \in \mathbb{R}^n$ is assumed to satisfy $S_c^+(v) \leq 2$, we have that $S_c^-(\Delta_c \sign(v)) \leq 2$ and $v \neq 0$. In order to prove our claim, we verify that \cref{item:PMP_ineq} of \cref{lem:VB_2}  is fulfilled for the following three subcases:
    \begin{enumerate}
        \item If $S^{+}_c(v) = 0$, it holds that $S_c^-(v) =0$. Hence, $\sign(v)$ can either be $\boldsymbol{1}_n$ or $-\boldsymbol{1}_n$. Therefore, $\Delta_c \sign(v) = 0$ satisfies \cref{item:PMP_ineq} of \cref{lem:VB_2} .
        \item If $S^{+}_c(v) = 2$ and $S^{-}_c(v) = 0$, then after a possible cyclic permutation, $\sign(v)$ can be either of the form $[\boldsymbol{1}_{n-1}^\top, 0]^\top$ or $[-\boldsymbol{1}_{n-1}^\top, 0]^\top$, in which case direct computations show that \cref{item:PMP_ineq} of \cref{lem:VB_2}  is fulfilled.
        \item If $S^{+}_c(v) = 2$ and $S^{-}_c(v) = 2$ we need to check all cases given in \eqref{eq:four_cases}. We will show this only for \eqref{eq:case_a} as the others follow similarly.
        Let $\sign(v) = [\boldsymbol{1}_{n_1}^\top, 0, -\boldsymbol{1}_{n_2}^\top, 0]^\top$ with $n_1+n_2+2 = n$ and $n_1, n_2 > 0$. Then, all possible minors of order $2$ of $H_{\Delta_c \sign(v)}$ are given by the submatrices
        \begin{align*}
            \begin{bmatrix}
                0 & 0 \\ 0 & 0
            \end{bmatrix},
            \pm\begin{bmatrix}
                0 & 1 \\ 0 & 0
            \end{bmatrix},
            \pm\begin{bmatrix}
                1 & 1 \\ 0 & 1
            \end{bmatrix},
            \pm\begin{bmatrix}
                0 & 0 \\ 1 & 0
            \end{bmatrix},
            \pm\begin{bmatrix}
                1 & 0 \\ 1 & 1
            \end{bmatrix}.
        \end{align*}
        As these minors are nonnegative and corrrespond to $(\Delta_c \sign(v_t))^2 - \Delta_c \sign(v_{t-1}) \Delta_c \sign(v_{t+1})$, we have verified \cref{item:PMP_ineq} of \cref{lem:VB_2}.
    \end{enumerate}
    Thus, in all three cases, \cref{lem:VB_2} applies, which proves our claim. \hfill $\Box$

\subsection{Proof of \cref{lem:Mickey_pres}}

    We begin by showing that $S_c^-(\Delta_c o^P) \leq 2$.
    Since $g$ is strictly monotonically decreasing on its support by \cref{ass:mono_dec_g}, we know that $S_c^-(\overline{g}^P)\leq S_c^+(\overline{g}^P) \leq 2$.
    Therefore, for any $u \in \ell_{\infty}(P)$ satisfying \cref{ass:osci_constr_1}, it holds by \cref{lem:sign_PM} that
    \begin{align*}
        S_c^-(\Delta_c o^P) & = S_c^-(-[\Delta \mathcal{C}_g(\sign(u))](0:P-1)) \\
        & = S_c^- \left( \Delta_c H_{\overline{g}^P}\sign(u^P) \right) \\
        & = S_c^- \left( \Delta_c H_{\sign(u^P)} \overline{g}^P \right) \leq 2.
    \end{align*}

    Next, we show that $S^{+}_c(o^P) \leq 2$ under the additional assumption that $P_p(u^P) = P_n(u^P)$. It suffices to verify the four cases given in \eqref{eq:four_cases}. We will do this only for \eqref{eq:case_a} using proof by contradiction, as the others follow similarly.

    Letting $\hat{u}^P := \sign(u^P) = [\boldsymbol{1}_{\frac{P}{2}-1}^\top, 0, -\boldsymbol{1}_{\frac{P}{2}-1}^\top, 0]^\top$ and suppose $S^-_c(\Delta_c o^P) \leq 2$ but $S^+_c(o^P) > 2$, then there exists $t_a \in (1:P-1)$ such that $o^P_{t_a} = o^P_{t_a +1} = 0$. Therefore,
    \begin{align}
        o^P_{t_a+1} - o^P_{t_a} & = -(H_{\overline{g}^P} \hat{u}^P)(t_a+1) + (H_{\overline{g}^P} \hat{u}^P)(t_a) \nonumber \\
        & = -(H_{\hat{u}^P} \overline{g}^P)(t_a+1) + (H_{\hat{u}^P} \overline{g}^P)(t_a) \nonumber \\
        & = -(\overline{g}^P)^\top
        \left( \begin{bmatrix}
            \hat{u}^P_{t_a+1} \\ \vdots \\ \hat{u}^P_{t_a - P+2}
        \end{bmatrix} - 
        \begin{bmatrix}
            \hat{u}^P_{t_a} \\ \vdots \\ \hat{u}^P_{t_a - P+1}
        \end{bmatrix} 
         \right) \nonumber \\
        & = -(\overline{g}^P)^\top (Q_P^{t_a+1} \hat{u}^P - Q_P^{t_a} \hat{u}^P) \nonumber \\
        & = -(\overline{g}^P)^\top Q_P^{t_a+1} v_0 = 0, \label{eq:minus_sum}
    \end{align}
    where $v_0 = \hat{u}^P - Q_P^{-1}\hat{u}^P =[1, 0, \dots, 0, -1, -1, 0, \dots,$ $ 0, 1]^\top$ and $Q_P$ is the cyclic shift back matrix defined in \eqref{eq:cyclic_back_shift_mat}. 
    Moreover, we also have that
    \begin{align}\label{eq:plus_sum}
         o^P_{t_a+1} + o^P_{t_a} & = -(\overline{g}^P)^\top (Q_P^{t_a+1} \hat{u}^P + Q_P^{t_a} \hat{u}^P) \nonumber \\
         & = -(\overline{g}^P)^\top Q_P^{t_a+1} v_1 = 0, 
    \end{align}
    where $v_1=\hat{u}^P + Q_P^{-1}\hat{u}^P =[1, 2, \cdots, 2, 1, -1, -2, \cdots,$ $ -2, -1]^\top$. Then, if $\hat{g}^P = (Q_P^{t_a+1})^{\top} \overline{g}^P$, there exists an $i \in (1:P)$ such that $\hat{g}^P_{i+1} > \dots > \hat{g}^P_P > \hat{g}^P_1 > \dots >  \hat{g}^P_{i}$ by \cref{ass:mono_dec_g} and, therefore,
    \begin{align*}
        & o^P_{t_a+1} - o^P_{t_a} > 0, \text{ if } 1 \leq i \leq P_p(u^P), \\
        & o^P_{t_a+1} + o^P_{t_a} > 0, \text{ if } i = P_p(u^P)+1, \\
        & o^P_{t_a+1} - o^P_{t_a} < 0, \text{ if } P_p(u^P)+2 \leq i \leq P-1, \\
        & o^P_{t_a+1} + o^P_{t_a} < 0, \text{ if } i = P.
    \end{align*}
    These conditions show that the equalities in \eqref{eq:minus_sum} and \eqref{eq:plus_sum} cannot hold simultaneously.
    Hence, we conclude that $S_c^+(o^P) \leq 2$ and $S^-_c(\Delta_c o^P) \leq 2$. \hfill $\Box$

\subsection{Proof of \cref{thm:osci_pres_open_loop}}

    Let $P_p(u^P)=P_p(o^P)$ and $P_n(u^P)=P_n(o^P)$ hold. In order to prove that $P_p(u^P) = P_n(u^P)$, we show that neither $P_p(u^P) > P_n(u^P)$ nor $P_p(u^P) < P_n(u^P)$ can arise.

    To this end, it suffices to consider the four cases in \eqref{eq:four_cases}. We only demonstrate this for \eqref{eq:case_a} using proof by contradiction, as the others follow similarly.

    Assume that $\sign(u^P) = [\boldsymbol{1}_{\frac{P}{2}-1}^\top, 0, -\boldsymbol{1}_{\frac{P}{2}-1}^\top, 0]^\top$ with $P_p(u^P) > P_n(u^P)$. By \cref{lem:Mickey_pres}, we have $S_c^-(\Delta_c u^P) \leq 2$ and $S^-_c(u^P) \leq 2$. Therefore, there exists a $\tau \in \mathbb{Z}$ such that
    \begin{align} \label{eq:time_delay_mul}
        \hat{o}^P := -Q_P^{\tau} H_{\overline{g}^P} \sign(u^P) = 
        \begin{bmatrix}
            \hat{o}^P_{+} \\ \hat{o}^P_\ominus
        \end{bmatrix},
    \end{align}
    where the circulant matrix $H_{\overline{g}^P} > 0$, $\hat{o}^P_{+} > 0$ and $\hat{o}^P_{\ominus} \leq 0$. 
    Note that, since the cyclic shift operater leaves $P_p$ and $P_n$ invariant, we have $P_p(\hat{o}^P) = P_p(Q_P^{\tau}o^P)$ and $P_n(\hat{o}^P) = P_n(Q_P^{\tau}o^P)$. Next, by \eqref{eq:conv_to_mat_mul}, we have $\hat{o}^P = -Q_P^{\tau} H_{\overline{g}^P} \sign(u^P) = -Q_P^{\tau} H_{\sign(u^P)} \overline{g}^P = [\hat{o}^P_1, \dots, \hat{o}^P_P]^\top$, where
    \begin{align} \label{eq:H-sign-u} 
        H_{\sign(u^P)} = \begin{bmatrix}
            1 & 0 & -1 & \cdots & -1 & 0 & 1 & \cdots & 1 \\
            1 & 1 & 0 & \cdots & -1 & -1 & 0 & \cdots & 1 \\
            \vdots & \vdots & \vdots &  & \vdots & \vdots & \vdots &   & \vdots \\
            -1 & -1 & -1 & \cdots & 1 & 1 & 1 &\cdots & 0 \\
            0 & -1 & -1 & \cdots & 0 & 1 & 1 &\cdots & 1 \\
        \end{bmatrix}.
    \end{align}
    Since $P_p(\hat{o}^P) = P_p(o^P) = P_p(u^P)$, it follows from $P_n(u^P)+1 \leq P_p(u^P)$ that $\hat{o}^P_{1} > 0$ and $\hat{o}^P_{P_n(u^P)+1}> 0$. However, since all elements in $\hat{g}^P$ are positive by $g$ being positive on its infinite support, it follows that
    \begin{align*}
        \hat{o}^P_{1} + \hat{o}^P_{P_n(u^P)+1} =  -2\hat{g}^P_1 - \hat{g}^P_2 - & \hat{g}^P_{2P_n(u^P)+3} \\& -\sum_{i = 2P_n(u^P)+4}^P 2\hat{g}^P_i < 0,
    \end{align*}
    which contradicts $P_n(u^P) = P_n(\hat{o}^P) = P_n(o^P)$. In a similar way, one can prove that $P_p(o^P) \neq P_p(u^P)$ when $P_p(u^P) < P_n(u^P)$.
    
    Therefore, if $P_p(u^P) = P_p(o^P)$ and $P_n(u^P) = P_n(o^P)$ hold, then $P_p(u^P) = P_n(u^P)$. Furthermore, if $P_z(u^P)=0$, then $P_z(o^P)=P_z(u^P)=0$ and $P = 2 P_p(u^P) = 2 P_p(o^P)$. \hfill $\Box$

\subsection{Proof of \cref{thm:non_oscillation}}

    From \cref{thm:osci_pres_open_loop}, the closed-loop system cannot admit a self-oscillation $u$ satisfying \cref{ass:osci_constr_1} if $P_p(u^P) \neq P_n(u^P)$.
    Hence, suppose that $u^P$ satisfies \cref{ass:osci_constr_1} with $P_p(u^P) = P_n(u^P)$ and $u^P = -H_{\overline{g}^P} \sign(u^P)$. We need to examine the four cases in \eqref{eq:four_cases}, which we will only do \eqref{eq:case_a}, as the others follow similarly.

    To this end, we begin by noticing that if $u$  in \eqref{eq:case_a} fulfills $u=-\mathcal{C}_g(\sign(u))$, it also need to hold that
    \begin{align*}
        u^P 
        = -H_{\overline{g}^P} \sign(u^P)
        = [{u^P_+}^\top, 0, {u^P_-}^\top, 0]^\top,
    \end{align*}
    where $\sign(u^P) = [\boldsymbol{1}_{\frac{P}{2}-1}^\top, 0, -\boldsymbol{1}_{\frac{P}{2}-1}^\top, 0]^\top$, $u^P_+ > 0$ and $u^P_- < 0$. Since we also have the relation $u^P=-H_{\sign(u^P)}\bar{g}^P$, where $H_{\sign(u^P)}$ is defined in \eqref{eq:H-sign-u}. 
    Then, since $g(0)>0$ and $g$ satisfies \cref{ass:mono_dec_g}, $\overline{g}^P_1 > \overline{g}^P_2 > \dots > \overline{g}^P_P$, which is why
    \begin{align*}
        0 < u^P_{P_p(u^P)} = -\sum_{i = 1}^{P_p(u^P)} \overline{g}^P_i + \sum_{i = P_p(u^P)+2}^{2P_p(u^P)+1} \overline{g}^P_i < 0.
    \end{align*} 
    This yields the desired contradiction. 
    \hfill $\Box$

\subsection{Proof of \cref{thm:bound_max_osc}}
    We need to verify our claim for the four main cases in \eqref{eq:four_cases}.
    Since by \cref{thm:osci_pres_open_loop}, $P_p(u^{P})=P_n(u^{P})$, it follows that $P$ satisfying $P \geq P_d$ must even in case of \eqref{eq:case_a} and \eqref{eq:case_d} and odd for \eqref{eq:case_b} and \eqref{eq:case_c}. In the following, we will consider the odd and even cases separately.
    
    We begin with the case, where $P = 2n \geq P_d$, $n \in \mathbb{Z}_{>0}$. In case of \eqref{eq:case_a}, we have that $\sign(u^{2n}) = [\boldsymbol{1}^{\top}_{n-1}, 0, -\boldsymbol{1}^{\top}_{n-1}, 0]^\top$ with $u^{2n} = -Q_{2n}^{P_d}y^{2n}$ and $y^{2n} = H_{\overline{g}_0^{2n}} \sign(u^{2n})$. Since $g_0$ satisfies \cref{ass:mono_dec_g}, the vector $\overline{g}_0^{2n} = [\overline{g}^{2n}_{0,1}, \dots, \overline{g}^{2n}_{0,2n}]^\top$ is strictly monotonically decreasing and, therefore, $y^{2n}_{n} = (\overline{g}^{2n}_{0, n} + \dots + \overline{g}^{2n}_{0, 2}) - (\overline{g}^{2n}_{0, 2n} + \dots + \overline{g}^{2n}_{0, n+2}) >0$. Thus, since $u^{2n} = -Q_{2n}^{P_d}y^{2n}$, this implies that $u^{2n}_{n+P_d} = -y^{2n}_{n} < 0$. However, as $u^{2n}_t < 0$ for $t \in (n+1:2n-1)$, it follows that $n+P_d \leq 2n-1$, or equivalently, $n \geq P_d +1$, which is why $P = 2n \geq 2P_d+2$. 

    For \eqref{eq:case_d}, we have that $\sign(u^{2n}) = [\boldsymbol{1}^\top_{n}, -\boldsymbol{1}^\top_{n}]^\top$. Then, analogues arguments yield that $u^{2n}_{n+P_d} = -y^{2n}_{n} < 0$ and, therefore, $n+P_d \leq 2n$. Hence, $P \geq 2P_d$ also in this case.
    
    Next, we consider the case of odd $P$ and since the proofs for \eqref{eq:case_b} and \eqref{eq:case_c} are similar, we restrict ourselves to \eqref{eq:case_b}. To this end, let $P = 2n+1 > P_d$ and $\sign(u^{2n+1}) = [\boldsymbol{1}_{n}^\top, 0, -\boldsymbol{1}_{n}^\top]^\top$ such that $u^{2n+1} = -Q_{2n+1}^{P_d}y^{2n+1}$ with $y^{2n+1} = H_{\overline{g}_0^{2n+1}} \sign(u^{2n+1})$. Since $P_d > 1$ and $g_0(0)>0$, the strict monotonicity of $\overline{g}_0^{2n+1}$ ensures $u^{2n+1}_{n+1+P_d} = -y^{2n+1}_{n+1} < 0$. Thus, $n+1+P_d \leq 2n+1$, i.e., $P \geq 2P_d+1$.
    
    Thus, in the all cases, we have shown that $P \geq 2P_d$ and that $P = 2P_d$ can only occur in case of \eqref{eq:case_d}.

    Next, we use proof by contradiction to show that $P \leq 2 (P_d + P_s)$. To this end, assume that $P \geq 2(P_d+P_s)+1$ with corresponding self-oscillation $u^{P} = [{u^{P}_+}^\top, {u^{P}_\ominus}^\top]^\top$, where the assumptions $ u^{P}_+> 0$ and $u^{P}_\ominus \leq 0$ cover all four cases in \eqref{eq:four_cases}. Then, by \cref{thm:osci_pres_open_loop}, $P_p(u^{P})$ satisfies $P_p(u^{P}) \geq \lceil \frac{P}{2} \rceil -1 \geq \lceil \frac{2(P_s+P_d)+1}{2} \rceil -1 = P_s+P_d$. Since $u^{P} = -Q_{P}^{P_d}y^{P}$ with $y^{P} = H_{\overline{g}_0^{P}} \sign(u^{P})$, we have $u^{P}_{P_p(u^{P})} = -y^{P}_{P_k}$, where $P_k := P_p(u^{P}) - P_d \geq P_s$. Moreover,
    \begin{align*}
        y^{P}_{P_k} = \sum_{i=1}^{P_k} \overline{g}^{P}_{0,i} + \sum_{i=P_k+1}^{P} \sign(u^{P}_i) \overline{g}^{P}_{0,P+P_k+1-i}. 
    \end{align*} 
    Since $u^{P}_i \leq 0$ for all $i \in (P_p(u^{P})+1:P)$, this implies $y^{P}_{P_k} \geq \sum_{i=1}^{P_k} \overline{g}^{P}_{0,i} - \sum_{i=P_k+1}^{P} \overline{g}^{P}_{0,i}$.
    Thus, since $\overline{g}^{P}_{0,i} = \sum_{k\in \mathbb{Z}} g_0(i-1+kP)$ with $g_0(t) \geq 0$ for all $t \in \mathbb{Z}$ and $g_0 \not \equiv 0$, we have 
    \begin{align*}
        y^{P}_{P_k}> \sum_{i=0}^{P_k-1} g_0(i) - \sum_{i=P_k}^{\infty} g_0(i)  \geq \sum_{i=0}^{P_s-1} g_0(i) - \sum_{i=P_s}^{\infty} g_0(i) > 0.
    \end{align*}
    This implies that $u^{P}_{P_p(u^{P})} = -y^{P}_{P_k} < 0$, which contradicts the assumption that $u^{P}_{P_p(u^{P})} > 0$. Hence, we must have $P \leq 2(P_d + P_s)$.

    Finally, we need to prove that $P = 2P_d$ for all $P_d \geq 1$ if and only if $\mathcal{C}_{g_0} (\sign{(u)}) (0) > 0$ with $\sign(u^{2P_d}) = [\boldsymbol{1}_{P_d}^\top, -\boldsymbol{1}_{P_d}^\top]^\top$. If $\mathcal{C}_{g_0} (\sign{(u)})(0) > 0$, then by strict monotonicity of $\overline{g}_0^{2P_d}$, we have 
    \begin{align*}
        0<\mathcal{C}_{g_0} (\sign{(u)})(0) \leq \dots \leq \mathcal{C}_{g_0} (\sign{(u)})(P_d-1).
    \end{align*}
    Thus, it follows that 
    \begin{align*}
        \sign(H_{\overline{g}^{2P_d}_0} \sign(u^{2P_d})) & =
        \sign \left( \begin{bmatrix}
            \mathcal{C}_{g_0} (\sign{(u)})(0) \\
            \vdots \\
            \mathcal{C}_{g_0} (\sign{(u)})(2P_d-1)
        \end{bmatrix} \right) \\
        & = [\boldsymbol{1}_{P_d}^\top,-\boldsymbol{1}_{P_d}^\top]^\top,
    \end{align*}
    so that $\sign \left(-Q_{2P_d}^{P_d} H_{\overline{g}^{2P_d}_0} \sign(u^{2P_d})\right) = [\boldsymbol{1}_{P_d}^\top,-\boldsymbol{1}_{P_d}^\top]^\top = \sign(u^{2P_d})$.
    Thus, $P = 2P_d$. Conversely, if $P = 2P_d$, then from $u^{2P_d} = -Q_{2P_d}^{P_d} H_{\overline{g}^{2P_d}_0} \sign(u^{2P_d})$ it follows that $\mathcal{C}_{g_0} (\sign{(u)}) (0) = u^{2P_d}_{P_d+1} < 0$.
    \hfill $\Box$

\subsection{Proof of \cref{cor:max_per_to_other}}
    Since $P = \frac{2P_d}{2n+1}$, $n \in \mathbb{Z}_{\geq 0}$, has to be even whenever $P \in \mathbb{Z}_{>0}$, it suffices to prove our claim for the cases of \eqref{eq:case_d}. By \cref{thm:bound_max_osc}, $\sign(u^{2P_d}) = [\boldsymbol{1}_{P_d}^\top, -\boldsymbol{1}_{P_d}^\top]^\top$ and 
    \begin{align*}
        \mathcal{C}_{g_0} (\sign(u))(0) = \overline{g}^{2P_d}_{0,1} - \sum_{i=2}^{P_d+1} \overline{g}^{2P_d}_{0,i} + \sum_{i=P_d+2}^{2P_d} \overline{g}^{2P_d}_{0,i} > 0.
    \end{align*} 
    Now, if $\sign(u^P) = [\boldsymbol{1}^\top_{\frac{P}{2}}, -\boldsymbol{1}^\top_{\frac{P}{2}}]^\top$, then by the monotonicity of $\overline{g}^{2P_d}_0$ and the relation 
    \begin{align*}
        \overline{g}^{P}_{0,i} & = \sum_{j \in \mathbb{Z}} g(i-1+jP) = \sum_{j \in \mathbb{Z}} g(i-1+j\frac{2P_d}{2n+1})\\
        & = \sum_{l \in \mathbb{Z}} \sum_{m=0}^{2n} g\left(i - 1 + ((2n+1)l+m)\frac{2P_d}{2n+1}\right) \\
        & = \sum_{m=0}^{2n} \sum_{l \in \mathbb{Z}} g\left(i+ \frac{2mP_d}{2n+1} -1+ 2l P_d\right) \\
        & = \sum_{m=0}^{2n} \sum_{l \in \mathbb{Z}} g\left(i + mP -1 + 2l P_d\right) = \sum_{m=0}^{2n} \overline{g}^{2P_d}_{0,i+mP},
    \end{align*} we obtain 
    \begin{align*}
        \mathcal{C}_{g_0} (\sign(u))(0) = & \overline{g}^{P}_{0,1} - \sum_{i=2}^{\frac{P}{2}+1} \overline{g}^{P}_{0,i} + \sum_{i=\frac{P}{2}+2}^{P} \overline{g}^{P}_{0,i} \\
        = &  \sum_{k=0}^{2n} \overline{g}^{2P_d}_{0,1+kP} - \sum_{i=2}^{\frac{P}{2}+1} \overline{g}^{2P_d}_{0,i+kP} \\
        & \quad \quad \quad \quad \quad + \sum_{i=\frac{P}{2}+2}^{P} \overline{g}^{2P_d}_{0,i+kP}  \\
        > & \overline{g}^{2P_d}_{0,1} - \sum_{i=2}^{P_d+1} \overline{g}^{2P_d}_{0,i} + \sum_{i=P_d+2}^{2P_d} \overline{g}^{2P_d}_{0,i} > 0.
    \end{align*}
    Since $\overline{g}_0^{2P_d}$ is also strictly monotonically decreasing, we further have $\mathcal{C}_{g_0} (\sign{(u)})(0) \leq \mathcal{C}_{g_0} (\sign{(u)})(1) \leq \dots \leq \mathcal{C}_{g_0} (\sign{(u)})(\frac{P}{2}-1)$. Thus, by \cref{thm:osci_pres_open_loop}
    \begin{align*}
        - Q_{P}^{\frac{P}{2}} H_{\overline{g}^{P}_0} \sign(u^P)  & = -Q_{P}^{\frac{P}{2}}\begin{bmatrix}
            \mathcal{C}_{g_0} (\sign(u))(0) \\
            \vdots \\
            \mathcal{C}_{g_0} (\sign(u))(P-1)
        \end{bmatrix} \\
        & = - Q_{P}^{\frac{P}{2}} \begin{bmatrix}
            \boldsymbol{1}_{\frac{P}{2}} \\ -\boldsymbol{1}_{\frac{P}{2}}
        \end{bmatrix}
        = \begin{bmatrix}
            u^P_+ \\ u^P_-
        \end{bmatrix},
    \end{align*}
    where $u^P_+ > 0$ and $u^P_- < 0$. Then, by $Q_{P}^{\frac{P}{2}} = Q_{P}^{\frac{(2n+1)P}{2}} = Q_{P}^{P_d}$, we deduce that $u^{P} = - Q_{P}^{\frac{P}{2}} H_{\overline{g}^{P}_0} \sign(u^P) = - Q_{P}^{P_d} H_{\overline{g}^{P}_0} \sign(u^P)$.
    Thus, the closed-loop system admits a $P$-periodic self-oscillation. \hfill $\Box$

\subsection{Proof of \cref{cor:convex_bound}}

From \cref{thm:bound_max_osc} we already know that $P \geq 2P_d$ when $P \geq P_d$. Next, we use proof by contradiction to show that $P \leq 4P_d +2$. 
It suffices to show this for the four cases in \eqref{eq:four_cases}. For even $P$, we consider the forms of \eqref{eq:case_a} and \eqref{eq:case_d} and for odd $P$ the forms of \eqref{eq:case_b} and \eqref{eq:case_c}.

To this end, assume that $P \geq 4P_d + 4$ is even and $u^{P} = -Q_{P}^{P_d} H_{\overline{g}_0^{P}} \sign(u^{P})$. We want to prove that $u^{P}_{\lceil{\frac{P}{4}} \rceil + P_d} \leq 0$ for the cases \eqref{eq:case_a} and \eqref{eq:case_d}.

For \eqref{eq:case_d}, let $\sign(u^{P}) = [\boldsymbol{1}^{\top}_{\frac{P}{2}},-\boldsymbol{1}^{\top}_{\frac{P}{2}}]^{\top}$.
Since $g_0$ satisfies \cref{ass:mono_dec_g}, we obtain 
\begin{align*}
    u^{P}_{\lceil{\frac{P}{4}} \rceil + P_d} = -\sum_{i=1}^{\lceil{\frac{P}{4}} \rceil} \overline{g}^{P}_{0,i} + \sum_{i=\lceil{\frac{P}{4}} \rceil +1}^{\lceil{\frac{P}{4}} \rceil + \frac{P}{2}} \overline{g}^{P}_{0,i} - \sum_{i=\lceil{\frac{P}{4}} \rceil + \frac{P}{2}+1}^{P} \overline{g}^{P}_{0,i},
\end{align*} 
and, thus,
\begin{align*}
    u^{P}_{\lceil{\frac{P}{4}} \rceil + P_d} \leq \sum_{i=1}^{\frac{P}{2} - \lceil{\frac{P}{4}} \rceil} -\overline{g}^{P}_{0, i} - \overline{g}^{P}_{0, P+1-i} + \overline{g}^{P}_{0, \frac{P}{2}+i} + \overline{g}^{P}_{0, \frac{P}{2}-i+1},
\end{align*}
where the equality holds exactly when $\frac{P}{4} = \lceil \frac{P}{4} \rceil$; otherwise, an extra item $-\overline{g}^{P}_{0,\lceil \frac{P}{4} \rceil} + \overline{g}^{P}_{0,\lceil \frac{P}{4} \rceil+\frac{P}{2}} < 0$ appears, ensuring the inequality is strict. 
By the convexity of $g_0$ on its support, for every $i \in (1:\frac{P}{2} - \lceil{\frac{P}{4}} \rceil)$ we have
\begin{align*}
    \overline{g}^{P}_{0, i} + \overline{g}^{P}_{0, P+1-i} \geq \overline{g}^{P}_{0, \frac{P}{2}+i} + \overline{g}^{P}_{0, \frac{P}{2}-i+1}
\end{align*}
and, thus, $u^{P}_{\lceil{\frac{P}{4}} \rceil + P_d} \leq 0$.

Similarly, for the case \eqref{eq:case_a}, let $\sign(u^{P}) = [\boldsymbol{1}^{\top}_{\frac{P}{2}-1}, 0, -\boldsymbol{1}^{\top}_{\frac{P}{2}-1}, 0]^\top$, then
\begin{align*}
    u^{P}_{\lceil{\frac{P}{4}} \rceil + P_d} = & -\sum_{i=1}^{\lceil{\frac{P}{4}} \rceil} \overline{g}^{P}_{0,i} + \sum_{i=\lceil{\frac{P}{4}} \rceil +1}^{\lceil{\frac{P}{4}} \rceil + \frac{P}{2}} \overline{g}^{P}_{0,i} - \sum_{i=\lceil{\frac{P}{4}} \rceil + \frac{P}{2}+1}^{P} \overline{g}^{P}_{0,i} \\
    & + \overline{g}^{P}_{0,P-1} - \overline{g}^{P}_{0,\frac{P}{2}-1} \\
    < & \sum_{i=1}^{\lceil{\frac{P}{4}} \rceil} \overline{g}^{P}_{0,i} + \sum_{i=\lceil{\frac{P}{4}} \rceil +1}^{\lceil{\frac{P}{4}} \rceil + \frac{P}{2}} \overline{g}^{P}_{0,i} - \sum_{i=\lceil{\frac{P}{4}} \rceil + \frac{P}{2}+1}^{P} \overline{g}^{P}_{0,i} \leq 0.
\end{align*}
However, for the cases \eqref{eq:case_a} and \eqref{eq:case_d}, since $P \geq 4$ and $P$ is even, it must hold that
\begin{align*}
    \lceil{\frac{P}{4}} \rceil + P_d \leq \lceil{\frac{P}{4}} \rceil + \frac{P-3}{4} \leq \frac{P+2}{4} + \frac{P-3}{4} \leq \frac{2P-1}{4}.
\end{align*}
Further, $\lceil{\frac{P}{4}} \rceil + P_d \leq \frac{P}{2} -1$ since $\lceil{\frac{P}{4}} \rceil + P_d$ is an integer. 
Thus, it follows that $u^{P}_{\lceil{\frac{P}{4}} \rceil + P_d} > 0$, which contradicts with $u^{P}_{\lceil{\frac{P}{4}} \rceil + P_d} \leq 0$. This proves our claim for the cases \eqref{eq:case_a} and \eqref{eq:case_d}.

Next, we consider the case of odd $P \geq 4P_d + 3$.
To this end, let $P \geq 8$ be odd and $u^{P} = -Q_P^{P_d} H_{\overline{g}_0^{P}} \sign(u^{P})$.
First, we need to prove that $u^{P}_{\lceil{\frac{P}{4}} \rceil + P_d} \leq 0$ for case \eqref{eq:case_b} and \eqref{eq:case_c}. Since the proofs for \eqref{eq:case_b} and \eqref{eq:case_c} are similar, we restrict ourselves to \eqref{eq:case_b}. Let $\sign(u^{P}) = [\boldsymbol{1}^{\top}_{\frac{P-1}{2}}, 0, -\boldsymbol{1}^{\top}_{\frac{P-1}{2}}]^\top$ such that $u^{P} = -Q_P^{P_d} H_{\overline{g}_0^{P}} \sign(u^{P})$. Thus, using a similar argument, we have
\begin{align*}
    u^{P}_{\lceil{\frac{P}{4}} \rceil + P_d} = -\sum_{i=1}^{\lceil{\frac{P}{4}} \rceil} \overline{g}^{P}_{0,i} + \sum_{i=\lceil{\frac{P}{4}} \rceil +1}^{\lceil{\frac{P}{4}} \rceil + \frac{P-1}{2}} \overline{g}^{P}_{0,i} - \sum_{i=\lceil{\frac{P}{4}} \rceil + \frac{P+3}{2}}^{P} \overline{g}^{P}_{0,i} \leq 0,
\end{align*} 
where the equality holds excatly when $\lceil \frac{P}{4} \rceil = \frac{P-1}{4}$. However, for the cases \eqref{eq:case_b} and \eqref{eq:case_c}, since $P \geq 8$ and $P$ is odd, it must hold that
\begin{align*}
    \lceil{\frac{P}{4}} \rceil + P_d \leq \lceil{\frac{P}{4}} \rceil + \frac{P-3}{4} \leq \frac{P+3}{4} + \frac{P-3}{4} \leq \frac{P}{2}.
\end{align*}
Further, $\lceil{\frac{P-1}{4}} \rceil + P_d \leq \frac{P-1}{2}$ since $\lceil{\frac{P-1}{4}} \rceil + P_d$ is an integer. Thus, it follows that $u^{P}_{\lceil{\frac{P-1}{4}} \rceil + P_d} > 0$.
Therefore, we have also shown our claim for the cases \eqref{eq:case_b} and \eqref{eq:case_c}.
\hfill $\Box$

\bibliographystyle{ieeetr}
\bibliography{ref.bib}

\end{document}